\documentclass[12pt,reqno]{amsart}

\usepackage{latexsym}
\usepackage{amssymb}

\renewcommand{\epsilon}{\varepsilon}

\newcommand{\var}{{\operatorname{Var}}}
\newcommand{\Var}{{\bf Var}}

\newcommand{\sm}{\smallsetminus}
\newcommand{\szego}{Szeg\H{o} }

\newcommand{\inv}{^{-1}}
\newcommand{\kahler}{K\"ahler }
\newcommand{\sqrtn}{\sqrt{N}}
\newcommand{\wt}{\widetilde}

\newcommand{\PP}{{\mathbb P}}

\newcommand{\R}{{\mathbb R}}
\newcommand{\C}{{\mathbb C}}

\newcommand{\Z}{{\mathbb Z}}

\newcommand{\CP}{\C\PP}
\renewcommand{\d}{\partial}
\newcommand{\dbar}{\bar\partial}
\newcommand{\ddbar}{\partial\dbar}

\newcommand{\E}{{\mathbf E}}

\newcommand{\half}{{\textstyle \frac 12}}
\newcommand{\vol}{{\operatorname{Vol}}}

\newcommand{\SU}{{\operatorname{SU}}}

\renewcommand{\phi}{\varphi}

\newcommand{\ccal}{\mathcal{C}}
\newcommand{\dcal}{\mathcal{D}}
\newcommand{\ecal}{\mathcal{E}}

\newcommand{\ical}{\mathcal{I}}

\newcommand{\ncal}{\mathcal{N}}
\newcommand{\ocal}{\mathcal{O}}

\newcommand{\al}{\alpha}

\newcommand{\ga}{\gamma}

\newcommand{\la}{\lambda}
\newcommand{\ep}{\varepsilon}
\newcommand{\de}{\delta}
\newcommand{\De}{\Delta}
\newcommand{\om}{\omega}
\newcommand{\Om}{\Omega}

\newtheorem{theo}{{\sc Theorem}}[section]

\newtheorem{cor}[theo]{{\sc Corollary}}

\newtheorem{lem}[theo]{{\sc Lemma}}
\newtheorem{prop}[theo]{{\sc Proposition}}

\newtheorem{defin}[theo]{{\sc Definition}}

\title[Number variance of random zeros on complex manifolds II]
{Number variance of random zeros on complex manifolds, II: smooth
statistics}

\author{Bernard Shiffman}
\author{Steve Zelditch}
\address{Department of Mathematics, Johns Hopkins University, Baltimore, MD
21218, USA} \email{shiffman@math.jhu.edu, zelditch@math.jhu.edu}

\thanks{Research of the first author partially supported by NSF grants
DMS-0100474 and DMS-0600982; research of the second author
partially supported by NSF grants  DMS-0302518 and DMS-0603850.}

\dedicatory{To Joseph J. Kohn on the occasion of his 75th birthday}

\begin{document}

\begin{abstract}
We consider the zero sets $Z_N$ of  systems of $m$ random polynomials of
degree $N$ in $m$ complex variables, and we give asymptotic formulas for
the random variables given by summing a smooth test function over $Z_N$.
Our asymptotic formulas show that the variances for these smooth
statistics  have the growth $N^{m - 2}$. We also prove analogues for the
integrals of smooth test forms over the subvarieties defined by $k < m$
random polynomials. Such  linear statistics of random zero sets are smooth
analogues of the random variables  given by counting the number of zeros
in an open set, which we proved elsewhere to have variances of order $N^{m
- 1/2}$. We use the variance asymptotics and off-diagonal estimates of
\szego kernels to extend an asymptotic normality result of Sodin-Tsirelson
to the case of  smooth linear statistics for zero sets of codimension one
in any dimension $m$.

\end{abstract}

\maketitle


\section{Introduction}

This article is concerned with zero sets of systems of Gaussian
random polynomials (or more generally,  of sections of a positive
holomorphic line bundle over a compact \kahler manifold $M_m$) as
the degree $N \to \infty$.  One of the most fundamental
statistical quantities  is the number
$\ncal_N^U(p_1^N,\dots,p_m^N)$  of  zeros in a bounded open set $U
\subset \C^m$ of a  system $\{p^N_1,\dots,p^N_m\}$ of $m$
independent Gaussian random polynomials. The expected value of
this random variable was shown in \cite{SZ} to be the integral of
the \kahler volume form over $U$ (times a universal constant).  In
a recent article \cite{SZa}, we gave an asymptotic formula for the
variance of this random variable. We also give analogous results
for the volume of the simultaneous zero set of $k<m$ polynomials
or sections. In this article  we apply the methods of \cite{SZa}
to the analogous `smooth linear statistics', i.e. the sum (or
integral) of a smooth test function over the zeros of a system of
random polynomials. Such smooth linear statistics arise as smooth
approximations for discontinuous random variables such as
$\ncal_N^U$  and also arise in a number of other problems (see the
discussion in \S \ref{OPEN}). Our main results give asymptotic
formulas for the variance of these  smooth linear statistics.
(Mean value asymptotics for these random variables were given in
\cite{SZ}.) As may be expected, the variances are of much lower
order in the degree $N$ than in the non-smooth case. Further, we
prove that in the codimension one case, the smooth linear
statistics are asymptotically normal, extending a result of
Sodin-Tsirelson \cite{ST}.

 To state our results precisely, we need some notation
and background.  We let $(L,h)\to M$ be a positively curved
Hermitian holomorphic line bundle over a compact complex manifold
$M$, and we give $H^0(M,L^N)$ the {\it Hermitian Gaussian
measure\/} induced by $h$ and the \kahler form $\om= \frac i2
\Theta_h$ (see Definition~\ref{HG}).  For  $m$ independent random
sections $s_j^N\in H^0(M,L^N)$, $1\le j\le m$, the number of
simultaneous zeros of the sections in a smooth domain $U \subset
M$ is given by
$$\ncal^U_N(s_1^N,\dots,s_m^N):=\#\{z\in U:s_1^N(z)= \cdots
=s_m^N(z)=0\}. $$ In \cite{SZ}, we proved that the expected value
$\E \big(\ncal^U_N\big)$ of the random variable $\ncal^U_N$ has the
asymptotics
\begin{equation} \label{MEAN} \E \big(\ncal^U_N\big) =
\frac {N^m}{\pi^m}\,\int_U \om^m+O(N^{m-1}), \end{equation} and in
\cite{SZa}, the variance of the random variable is shown to have
the asymptotics,
\begin{equation}\label{number}\var\big(\ncal^U_N\big) =
N^{m-1/2} \,\left[\nu_{mm} \,\vol_{2m-1}(\d U)
 +O(N^{-\frac 12 +\ep})\right]\;,\end{equation} where $\nu_{mm}$ is a
universal positive constant (\cite[Theorem~1.1]{SZa}). Analogous
results are proved for the volumes of zero sets  of  $k\le
m-1$ independent random sections $s_j^N\in H^0(M,L^N)$, $1\le j\le
k$: in this case,
\begin{equation}\label{volume}\var\big(\vol_{2m-2k}[Z_{s_1^N,\dots,s_k^N}\cap
U]\big) = N^{2k-m-1/2} \,\left[\nu_{mk} \,\vol_{2m-1}(\d U)
 +O(N^{-\frac 12 +\ep})\right]\;,\end{equation} where $\nu_{mk}$ is a
universal constant (\cite[Theorem~1.4]{SZa}); in particular,
$\nu_{m1} = \frac{\pi^{m-5/2}}{8}\;\zeta(m+\textstyle \half)$.
More generally, the domain could be piecewise smooth without cusps.

In this article, we are interested in the smooth analogue of
$\ncal^U_N$ where we integrate a smooth test function rather than
the characteristic function of a smooth domain over the zero set.
Given a test function $\phi \in \dcal(M)$, we consider the random
variable
\begin{equation*}\label{SLS} (Z_{s_1^N,\dots,s_m^N},\phi)= \sum_{
s_1^N(z)=\cdots = s_m^N(z)=0} \phi(z).  \end{equation*} When the
system is not full, we  define
\begin{equation*}\label{k} (Z_{s_1^N,\dots,s_k^N},\phi)= \int_{
s_1^N(z)=\cdots = s_k^N(z)=0} \phi(z), \qquad \phi\in
\dcal^{m-k,m-k}(M)\;. \end{equation*} The expected value of
$(Z_{s_1^N,\dots,s_m^N},\phi)$ is given by (see \eqref{indepE})
\begin{equation}\label{ave}\E(Z_{s_1^N,\dots,s_k^N},\phi)\approx
N^k\,\pi^{-k}\int_M \om^k\wedge\phi. \end{equation}The main result
of this article is an asymptotic formula for the variance:

\begin{theo}\label{sharp} Let $(L,h)$ be a
positive Hermitian holomorphic line bundle over a compact \kahler
manifold $(M,\om)$, where $\om= \frac i2 \Theta_h$, and let $1\le
k\le m$.  We give $H^0(M,L^N)$ the Hermitian Gaussian measure
induced by $h,\om$ (see Definition \ref{HG}).

Let  $\phi$ be a real  $(m-k,m-k)$-form on $M$ with $\ccal^3$ coefficients.
Then for independent random sections
$s_1^N,\dots,s_k^N\in H^0(M,L^N)$, we have
$$\var\big(Z_{s_1^N,\dots,s_k^N},\phi\big) = N^{2k-m-2}
\left[
\int_M B_{mk}\left(\ddbar\phi,\ddbar\phi\right)\Om_M
 +O(N^{-\frac 12 +\ep})\right]\;,$$  where $\Om_M$ is the volume form on $M$,
and $B_{mk}$ is a universal Hermitian form on the bundle
$T^{*m-k+1,m-k+1}(M)$. When $k=1$, we have
$B_{m1}\left(f\Om_M,f\Om_M\right)=
 \frac{\pi^{m-2}\,\zeta(m+2)}{4}\,
|f|^2,$ and hence
$$\textstyle\var\big(Z_{s^N},\phi\big) = N^{-m}
\left[ \frac{\pi^{m-2}\,\zeta(m+2)}{4}\,\|\ddbar\phi\|_{L^2}^2
 +O(N^{-\frac 12 +\ep})\right]\;.$$
\end{theo}
 In particular, for
the complex curve case $m=1$, we note that $|\ddbar\phi| = \half
|\Delta \phi|$, and thus
\begin{equation}\var\big(Z_{s^N},\phi\big) = N\inv
\left[\frac{\zeta(3)}{16\pi}\,\|\Delta\phi\|_2^2
 +O(N^{-\frac 12 +\ep})\right]\;.\label{smooth1}\end{equation}
The leading term in \eqref{smooth1} was obtained by Sodin and Tsirelson \cite{ST} for the case of random
polynomials $s^N\in H^0(\CP^1,\ocal(N))$ as well as for random holomorphic functions on $\C$ and on the disk.
(The constant  $\frac{\zeta(3)}{16\pi}$  was given in a private communication from M. Sodin.)

Here we say that  $B_{mk}$ is universal if there exists a
Hermitian inner product $B^0_{mk}$ on\break $
T^{*m-k+1,m-k+1}_0(\C^m)$, independent of $M$ and $L$, such that
for all $w\in M$ and all unitary transformations
$\tau:T^*_0(\C^m)\to T^*_w(M)$, we have
$B_{mk}(w)=\tau_*B^0_{mk}$. The global inner product
$\big(\phi,\psi\big) = \int_M
B_{mk}\left(\ddbar\phi,\ddbar\psi\right)\Om_M$ is certainly
positive semi-definite on $\dcal^{m-k,m-k}(M)$, since the variance
is nonnegative. We believe  that, in fact, $B^0_{mk}$ is positive
definite on $T_0^{*m-k+1,m-k+1}(\C^m)$.  This follows for $k=1$
from the above formula for $B_{m1}$; one should be able to verify
positivity for $k>1$ by using  the expansion \eqref{many} in the
proof of Theorem \ref{sharp} to compute a precise formula for
$B^0_{mk}$.

Thus the variance of the `smooth statistic'
$(Z_{s_1^N,\dots,s_k^N},\phi)$ is of lower order than the variance
of the number and volume statistics given by
\eqref{number}--\eqref{volume},  as expected.  In view of
(\ref{ave}), it is also self-averaging in the sense  that its
fluctuations  are of smaller order than its typical values.

An application of our methods is an extension of the Sodin-Tsirelson
\cite{ST} asymptotic normality result for smooth
statistics  to general one-dimensional ensembles and to
codimension one zero sets in higher dimensions:
\begin{theo} \label{AN} Let $(L,h)\to(M,\om)$ be as
in Theorem \ref{sharp} and give $H^0(M,L^N)$ the Hermitian Gaussian measure
induced by $h,\om$. Let $\phi$ be a real  $(m-1,m-1)$-form on $M$ with
$\ccal^3$ coefficients, such that
 $\ddbar\phi\not\equiv 0$.  Then for random sections $s^N$ in
$H^0(M,L^N)$, the distributions of the random variables
$$ \frac{(Z_{s^N},\phi)-\E(Z_{s^N},\phi)}{\sqrt{\var(Z_{s^N},\phi)}}$$
converge weakly to the
standard Gaussian distribution
$\ncal(0, 1)$ as $ N \to \infty$.\end{theo}

Sodin and Tsirelson \cite{ST} obtained the asymptotics of Theorem
\ref{AN} for random
functions on $\C$, $\CP^1$, and the disk. The proof of Theorem
\ref{AN} is a relatively straightforward application of the fundamental
\szego kernel asymptotics underlying Theorem \ref{sharp} to the argument in
\cite{ST}. (One easily sees that the random variable $(Z_{s^N},\phi)$ is
constant for all $N$ if
$\ddbar\phi\equiv 0$.)

Substituting the values of the expectation and variance of
$(Z_{s^N},\phi)$ from \eqref{ave} and Theorem \ref{sharp},
respectively, we have:
\begin{cor} With the same notation and hypotheses as
in Theorem \ref{AN}, the distributions of the random variables
$
N^{m/2}(Z_{s^N}-\frac N\pi\,\om,\phi)$ converge weakly to
$\ncal(0, \sqrt{\kappa_m}\, \|\ddbar\phi\|_2)$ as $ N \to \infty$, where
$\kappa_m= \frac{\pi^{m-2}\,\zeta(m+2)}{4}$.
\end{cor}

Here, $\ncal(0, \sigma)$ denotes the (real)
Gaussian distribution of mean zero and variance
$\sigma^2$.

We now summarize the key ideas in the proofs in \cite{SZa} and in
this paper. The variance in Theorem \ref{sharp}, as well as the
number and volume variances in \cite{SZa}, can be expressed in
terms of  the variance currents $\Var(Z_{s_1^N,\dots,s_k^N})$ of
the random currents $Z_{s_1^N,\dots,s_k^N}$.  In joint work with
P.  Bleher in 2000 \cite{BSZ1}, we introduced a bipotential $Q_N$
for the  `pair correlation function' $K^N_{21}$ of the volume
density of zeros of random sections in $H^0(M,L^N)$; this
bipotential satisfies
\begin{equation*}  \Delta_z \Delta_w  Q_N (z,w) =
K_{21}^N(z,w)\,. \end{equation*} The bipotential $Q_N$ is a universal function of the
normalized  \szego kernel (see \eqref{PN} and \eqref{QN}). Sodin and
Tsirelson \cite{ST} obtained a variance formula as well as asymptotic
normality for zeros of certain model one-dimensional  random holomorphic
functions by implicitly using this bipotential.

In \cite{SZa}, we showed that $Q_N$ is actually a `pluri-bipotential' for
the codimension-one variance current; i.e.,
\begin{equation}\label{pluri} (i\ddbar)_z\,(i\ddbar)_w \,Q_N(z,w) = {\bf
Var}\big(Z_{s^N}\big)\;.\end{equation}  We further found a formula
(Theorem \ref{variant}) for the higher codimension variance current
$\Var(Z_{s_1^N,\dots,s_k^N})$ in terms of $Q_N$ and its derivatives of
order $\le 4$. We then applied the  off-diagonal asymptotics of the
\szego kernel $\Pi_N(z,w)$  in \cite{SZ2} to obtain asymptotics of the
bipotential
$Q_N(z,w)$ and then of the number variance \eqref{number} as well as the
volume variance \eqref{volume}.

In this paper, we begin by reviewing basic facts about the \szego kernel
and summarizing the asymptotics from \cite{SZa} of the bipotential $Q_N$
as $N\to \infty$ as well as $d(z,w)\to 0$. To illustrate our ideas, we
apply these asymptotics to \eqref{pluri} to derive the codimension one
formula (i.e., the case $k=1$) of Theorem \ref{sharp}. We then  prove in
\S \ref{explicit} a slight modification (Corollary \ref{BIPOTk}) of the
formula for the higher codimension variance, which we use in \S
\ref{higher} to prove Theorem \ref{sharp}.  In \S \ref{s-normality} we
apply our \szego kernel asymptotics  to prove Theorem \ref{AN}.
Finally,  we formulate some related open problems in \S \ref{OPEN}.

\section{Background}\label{background}

In this section we summarize results from \cite{SZa} used in this
paper.

 We let  $(L,h)$ be a Hermitian holomorphic
line bundle over a compact \kahler manifold $M$. We
consider a local holomorphic frame $e_L$ over a trivializing chart
$U$. If  $s = f e_L$ is a section of $L$ over $U$,  its Hermitian norm is
given by
$\|s(z)\|_h = a(z)^{-\frac 12}|f(z)|$ where \begin{equation}
\label{a} a(z) = \|e_L(z)\|_h^{-2}\;. \end{equation} The curvature form of
$(L,h)$ is given locally by
$$\Theta_h= \ddbar \log a\;,$$ and the
{\it Chern form\/} $c_1(L,h)$ is given by
\begin{equation}\label{chern}c_1(L,h)=\frac{\sqrt{-1}}{2 \pi}
\Theta_h=\frac{\sqrt{-1}}{2 \pi}\d\dbar\log a\;.\end{equation} The
current of integration $Z_s$ over the zeros of a section $s\in H^0(M,L)$
is then given by the {\it Poincar\'e-Lelong formula\/},
\begin{equation} Z_{s} =
\frac{\sqrt{-1}}{ \pi } \partial \bar{\partial}\log |f| = \frac{\sqrt{-1}}{
\pi } \partial
\bar{\partial}\log\left\|s\right\|_{h} + c_1(L,h)\;,
\label{Zs}
\end{equation} where the second equality is a consequence of
\eqref{a}--\eqref{chern}.

We now assume that the Hermitian metric $h$ has strictly positive
curvature and we give
$M$ the \kahler form \begin{equation}\label{omega}\om= \frac i2\Theta_h=\pi
c_1(L,h)\;.\end{equation}
Next we describe the  natural Gaussian probability measures on the spaces
$H^0(M,L^N)$ of holomorphic sections of tensor powers
$L^N=L^{\otimes N}$ of the line bundle $L$:

\begin{defin}\label{HG} Let $(L,h)\to (M,\om)$ be as above,
and let $h^N$ denote the Hermitian metric on $L^N$ induced by $h$.
We give $H^0(M,L^N)$ the inner product induced by the \kahler form
$\om$ and the Hermitian metric $h^N$:
\begin{equation}\label{inner}\langle s_1, \bar s_2 \rangle = \int_M h^N(s_1,
s_2)\,\frac 1{m!}\om^m\;,\qquad s_1, s_2 \in
H^0(M,L^N)\,\;.\end{equation}  The \ {\em Hermitian Gaussian measure}
on
$H^0(M,L^N)$ is the complex Gaussian probability measure $\ga_N$
induced by the inner product \eqref{inner}:
\begin{equation*}d\ga_N(s)=\frac{1}{\pi^m}e^
{-|c|^2}dc\,,\qquad s=\sum_{j=1}^{d_N}c_jS^N_j\,,\end{equation*}
where $\{S_1^N,\dots,S_{d_N}^N\}$ is an orthonormal basis for
$H^0(M,L^N)$. It is of course independent
of the choice of orthonormal  basis.
\end{defin}

 The Gaussian ensembles $(H^0(M,L^N),\ga_N)$  were also  studied in
\cite{SZ,SZ2,BSZ1,BSZ2}; for the case of polynomials in one variable, they
become the $\SU(2)$ ensembles in \cite{BBL, Han, NV, Zh}; for polynomials
in $m$ complex variables, they are the $\SU(m+1)$ ensembles (see, e.g.,
\cite{SZ,  Zr}).

We consider the diagonal \szego kernels
$$\Pi_N(z,z):=\sum_{j=1}^{d_N} \|
S^N_j(z)\|^2_{h^N(z)}\;,$$ where the $S_j^N$ are as in the above
definition.
It follows from the leading terms of the  asymptotic
expansion of the diagonal \szego kernel   of
\cite{Cat,Ti,Z} that
\begin{equation}\label{CZ}\Pi_N(z,z)= \frac
1{\pi^m}N^m(1+O(N\inv))\,.\end{equation}
The expected value of the zero current of a random
holomorphic section in $H^0(M,L^N)$ is given by  the  basic
formula
\begin{equation}\label{EZ}\E Z_{s^N} = \frac i{2\pi}\ddbar
\log \Pi_{N}(z,z) + \frac N\pi \omega\,,\end{equation} and the expected
values of simultaneous zero currents are given by
\begin{equation} \label{indepE}\E\big(Z_{s_1^N,\dots,s_k^N}\big) =
\left[\E\big(Z_{s^N})\right]^{\wedge k} = \left(\frac i{2\pi}\ddbar
\log \Pi_{N}(z,z) + \frac N\pi \omega\right)^k = \frac
{N^k}{\pi^k}\,\om^k+O(N^{k-1}) \;,\end{equation}
for $1\le k\le m$ (see \cite{SZ,SZa}).
The final equality of  \eqref{indepE} is a consequence of
the asymptotic formula \eqref{CZ}.

For the variance asymptotics, we need the properties of the off-diagonal
\szego kernel:
\begin{equation}\label{zw}
|\Pi_N(z,w)|:=  \left\|\sum_{j=1}^{d_N}
S^N_j(z)\otimes\overline{ S^N_j(w)}\right\|_{h^N(z)\otimes
h^N(w)}\;.\end{equation} In particular, our variance formulas are
expressed in terms of the {\it normalized \szego kernel}
\begin{equation}\label{PN} P_N(z,w):=
\frac{|\Pi_N(z,w)|}{\Pi_N(z,z)^\frac 12 \Pi_N(w,w)^\frac
12}\;.\end{equation}

In \cite{SZa}, we
used the off-diagonal asymptotics for $\Pi_N$ from \cite{SZ2} to
provide the off-diagonal estimates given below for the  normalized \szego
kernel
$P_N$. As in   \cite{BSZ1, SZ, SZ2, SZa}), we
obtained these asymptotics by identifying the line bundle \szego
kernel $|\Pi_N(z,w)|$ of \eqref{zw} with the absolute value of a scalar
\szego kernel
$\Pi_N(x, y)$ on the unit circle bundle $X \subset L^{-1} \to M$ associated
to the Hermitian metric $h$.

Our estimates are of
two types: (1) `near-diagonal' asymptotics (Proposition
\ref{better}) for
$P_N(z,w)$ where the distance
$d(z,w)$ between
$z$ and
$w$ satisfies an upper bound $ d(z,w)\le b\left(\frac
{\log N}{N}\right)^{1/2}$ ($b\in\R^+$); (2) `far-off-diagonal' asymptotics
(Proposition \ref{DPdecay}) where  $ d(z,w)\ge b\left(\frac
{\log N}{N}\right)^{1/2}$:

\begin{prop}\label{DPdecay} {\rm \cite[Prop.~2.6]{SZa}} Let $(L,h)\to
(M,\om)$ be as in Theorem
\ref{sharp}, and let $P_N(z,w)$ be the normalized \szego kernel for
$H^0(M,L^N)$ given by \eqref{PN}. For
$b>\sqrt{j+2k}$,
$j,k\ge 0$, we have
$$ \nabla^j
P_N(z,w)=O(N^{-k})\qquad \mbox{uniformly for }\ d(z,w)\ge
b\,\sqrt{\frac {\log N}{N}} \;.$$
\end{prop}
\medskip

Here, $\nabla^j$ stands for the $j$-th covariant derivative.
The normalized \szego kernel $P_N$ also satisfies Gaussian
decay estimates valid very close to the diagonal. To give this estimate,
we consider a local normal coordinate chart
$\rho:U,z_0\to\C^m,0$ centered at a point $z_0\in U\subset M$,  and we
write, by abuse of notation,
\begin{equation*}P_N(z_0+u,z_0+v):=
P_N(\rho\inv(u),\rho\inv(v))\;.\end{equation*}

\begin{prop} \label{better}  {\rm \cite[Prop.~2.7--2.8]{SZa}} Let
$P_N(z,w)$ be as in Proposition
\ref{DPdecay}, and let $ z_0\in M$.  For $\ep,b>0$, there are constants
$C_j=C_j({M,\ep,b})$, $j \ge 2$, independent of the point $z_0$, such that
$$\textstyle  P_N\left(z_0+\frac u{\sqrtn},z_0+\frac v{\sqrtn}\right) =
e^{-\frac 12 |u-v|^2}[1 + R_N(u,v)]\;,$$ where
$$ \begin{array}{c}|R_N(u,v)|\le \frac {C_2}2\,|u-v|^2N^{-1/2+\ep}\,, \quad
|\nabla R_N(u)|
\le C_2\,|u-v|\,N^{-1/2+\ep}\,,
\\[8pt] |\nabla^jR_N(u,v)|\le C_j\,N^{-1/2+\ep}\quad j\ge 2\,,\end{array}$$
for $|u|+|v|<b\sqrt{\log N}$.\end{prop}

\subsection{The pluri-bipotential for the variance}\label{bipot}

For
random codimension
$k$ zeros, we have the {\it variance current\/} of
$Z_{s_1^N,\dots,s_k^N}$:
\begin{equation}\label{vc}
\Var\big(Z_{s_1^N,\dots,s_k^N}\big) =  \E\big(Z_{s_1^N,\dots,s_k^N}
\boxtimes Z_{s_1^N,\dots,s_k^N}
\big) -
\E\big(Z_{s_1^N,\dots,s_k^N}\big)\boxtimes \E
\big(Z_{s_1^N,\dots,s_k^N}\big)\in \dcal'^{2k,2k}(M\times M).
\end{equation}
The variance for the
`smooth zero statistics' is given by:
\begin{equation}\label{var}\var\big(Z_{s_1^N,\dots,s_k^N},\phi\big)
=\left(\Var\big(Z_{s_1^N,\dots,s_k^N}\big),\;\phi\boxtimes
\phi\right)\;.\end{equation}
Here we  write
$$R\boxtimes S = \pi_1^*R \wedge \pi_2^*S \in \dcal'^{p+q}(M\times
M)\;, \qquad \mbox{for }\ R\in \dcal'^p(M),\ S\in \dcal'^q(M)\;,$$
where $\pi_1,\pi_2:M\times M\to M$ are the projections to the
first and second factors, respectively.

For a current $T$ on $M\times
M$, we shall write
$$\d T=\d_1 T+\d_2 T\,,\qquad \d_1 = \sum dz_j\; \frac {\d}{\d z_j }\;,
\quad \d_2 = \sum dw_j\; \frac {\d}{\d w_j }\;,$$ where
$z_1,\dots,z_m$ are local coordinates on the first factor, and
$w_1,\dots, w_m$ are local coordinates on the second factor of
$M\times M$. In particular, $\d_1(R\boxtimes S)= (\d R)\boxtimes S$ and
$\d_2(R\boxtimes S)= R\boxtimes (\d S)$. We similarly
write
$$ \dbar T = \dbar_1 T +\dbar_2T\;.$$

In \cite{SZa}, we constructed a {\it pluri-bipotential\/} for
the variance current in codimension one, i.e.\  a function $Q_N\in
L^1(M\times M)$ such that
\begin{equation}\label{varcur}{\bf Var}\big(Z_{s^N}\big)=
-\d_1\dbar_1\d_2\dbar_2Q_N= (i\ddbar)_z\,(i\ddbar)_w \,Q_N(z,w)
\;.\end{equation} To describe our pluri-bipotential $Q_N(z,w)$, we define the
function
\begin{equation}\label{Gtilde} \wt G(t):= -\frac 1{4\pi^2}
\int_0^{t^2} \frac{\log(1-s)}{s}\,ds\ =\ \frac 1{4\pi^2}
\sum_{n=1}^\infty\frac{t^{2n}}{n^2}\;,
\qquad 0\le t\le 1.\end{equation}  Alternatively,
\begin{equation}\label{Gtilde1}\wt G(e^{-\la}) = -\frac
1{2\pi^2} \int_\la^\infty \log(1-e^{-2s})\,ds\;,\qquad
\la\ge 0\;.\end{equation}
\begin{theo} \label{BIPOT} {\rm \cite[Theorem 3.1]{SZa}} Let $(L,h)\to(M,\om)$ be as in
Theorem
\ref{sharp}. Let $Q_N:M\times M\to [0,+\infty)$ be
the function given  by
\begin{equation}
\label{QN} Q_N(z,w)= \wt G(P_N(z,w)) = -\frac 1{4\pi^2}
\int_0^{P_N(z,w)^2} \frac{\log(1-s)}{s}\,ds\;, \end{equation}
where $P_N(z,w)$ is the normalized \szego kernel given by
\eqref{PN}. Then $${\bf Var}\big(Z_{s^N}\big)=
 -\d_1\dbar_1\d_2\dbar_2Q_N\;.$$
\end{theo}
Theorem \ref{BIPOT} says that
\begin{equation} \label{varint1} \var( Z_{s^N},\phi) =
\big(-\d_1\dbar_1\d_2\dbar_2Q_N,\;\phi\boxtimes \phi\big)= \int_{M\times
M} Q_N(z,w)\,i\ddbar \phi(z)\wedge i\ddbar \phi(w)
\;,\end{equation} for all real $(m-1,m-1)$-forms $\phi$ on $M$
with $\ccal^2$ coefficients.

Since $P_N^2\in\ccal^\infty(M\times M)$ and $P_N(z,w)<1$ for $z\neq w$,
for sufficiently large $N$ (so that the Kodaira map for $L^N$ is an
embedding), it follows from \eqref{Gtilde} that  $Q_N$ is $\ccal^\infty$
off the diagonal, for $N \gg 0$.

Proposition \ref{DPdecay} implies that the pluri-bipotential decays rapidly
away from the diagonal:

\begin{lem}\label{Qdecay} {\rm \cite[Lemma 3.4]{SZa}} For $b>\sqrt{j+q+1},\ j\ge 0$,
we have
$$|\nabla^j Q_N( z,w)| =O\left(\frac
1{N^{q}} \right)\;,\quad \mbox{for }\ d(z,w)>\frac{b\sqrt{\log
N}}{\sqrtn}\;.$$ \end{lem}

Proposition \ref{better} yields the near-diagonal asymptotics:
\begin{lem}\label{FNnear} {\rm \cite[Lemma 3.5]{SZa}} For $b\in \R^+$, we have $$
Q_N\Big(z_0,z_0+\frac v {\sqrtn}\Big)= \wt G(e^{-\frac 12 |v|^2})+
O(N^{-1/2+\epsilon})\;,\qquad \mbox {for }\ |v|\le b\sqrt{\log N
}.$$\end{lem}

Recalling \eqref{Gtilde1}, we write,
\begin{equation}\label{F} F(\la):=\wt G(e^{-\la}) = -\frac
1{2\pi^2} \int_\la^\infty \log(1-e^{-2s})\,ds \qquad\quad
(\la\ge 0)\;,\end{equation} so that
$Q_N=F\circ(- \log P_N)$.  By Proposition \ref{better},
\begin{equation}\label{RN0}-\log P_N\left(z_0,z_0+\frac v{\sqrtn}\right)
={\half |v|^2} + O(|v|^2N^{-1/2+\ep})\quad \mbox{for }\
|v|<b\sqrt{\log N}\;.\end{equation} It follows from Lemma \ref{FNnear} and
\eqref{RN0} that $Q_N\in \ccal^1(M\times M)$ and the
first partial derivatives of $Q_N$ vanish along the diagonal in $M\times
M$, for $N\gg 0$.  (We note that $Q_N$ is $\ccal^\infty$ off the diagonal,
but is not
$\ccal^2$ at all points on the diagonal in
$M\times M$, as the computations in \cite{SZa} show.) We furthermore have
the near-diagonal asymptotics:

\begin{lem} \label{DISTM2} {\rm \cite[Lemma 3.7]{SZa}} There exist a
constant
$C_m\in\R^+$ (depending only on the dimension $m$) and $N_0=N_0(M)\in\Z^+$
such that for
$N\ge N_0$, we have:
\begin{itemize} \item[i)] The coefficients of the current
$\dbar_1\dbar_2 Q_N $ are locally bounded functions (given by
pointwise differentiation of $Q_N$), and we have the pointwise estimate
 $$|\dbar_1\dbar_2 Q_N(z,w) |\le C_mN\quad \mbox {for }\ 0<|w-z|<
b\sqrt{\frac {\log
 N}{N}}.
$$

\item[ii)] If $m\ge 2$, the  coefficients of the current
$\d_1\dbar_1\d_2\dbar_2Q_N$ are locally $L^{m-1}$ functions, and we
have the estimate
 $$|\d_1\dbar_1\d_2\dbar_2Q_N(z,w)| \le \frac{C_mN}{|w-z|^2}
\quad \mbox {for }\ 0<|w-z|< b\sqrt{\frac {\log N}{N}}.
$$
\end{itemize}\end{lem}

\begin{lem}\label{d4Qas} {\rm \cite[Lemma 3.9]{SZa}} For $N$ sufficiently large,
\begin{equation}\label{d4Qv} \textstyle-
\d_1\dbar_1\d_2\dbar_2Q_N(z_0,z_0+ \frac v{\sqrtn}) =
N\,\Var_\infty^{z_0}(v) +O\left(|v|^{-2}\,N^{1/2+\ep}\right)
\quad\mbox{for }\ 0<|v|<b\sqrt{\log N} \;,\end{equation} where
$\Var_\infty^{z_0}\in T^{*1,1}_{z_0}(M)\otimes \dcal'^{1,1}(\C^m)$
is given by
\begin{eqnarray}\label{many} \textstyle\Var_\infty^{z_0}(v) &:=&
\textstyle -\frac {1}{16}F^{(4)}(\half |v|^2)\ (\bar v\cdot
dz)(v\cdot d\bar z)(\bar v\cdot dv)(v\cdot d\bar v)\nonumber
\\&&-\ \textstyle \frac{1}{8} F^{(3)}(\half |v|^2)\, \big[(dz\cdot
d\bar z)(\bar v\cdot dv)(v\cdot d\bar v) +(v\cdot d\bar z)(\bar
v\cdot dv)(dz\cdot d\bar v)\nonumber\\&&\textstyle \qquad\qquad +
(\bar v\cdot dz) (d\bar z\cdot dv)( v\cdot d\bar v) + (\bar v\cdot
dz)(v\cdot d\bar z)(dv\cdot d\bar v) \big]\nonumber\\&&-\
\textstyle \frac 14 \,F''(\half |v|^2)\,\big[ (d\bar z\cdot
dv)(dz\cdot d\bar v) + (dz\cdot d\bar z)(dv\cdot d\bar
v)\big]\,.\end{eqnarray}
\end{lem}

Differentiating \eqref{F}, we note
that
\begin{equation}\label{F''} F''(\la)= \frac 1{\pi^2} \,\frac
1{e^{2\la}-1}\,,\quad F^{(3)}(\la)= -\frac 1{2\pi^2}\,
\mbox{csch}^2\la\,,\quad F^{(4)}(\la)= \frac 1{\pi^2}\,\coth
\la\,\mbox{csch}^2\la
\;.\end{equation}
Thus,  $$F^{(j)}(\la) = (-1)^j\, \frac{(j-2)!}{2\pi^2}\,\la^{-j+1} +O(1)
\qquad (\la>0)\,,$$ for $j\ge 2$, and hence
\begin{equation}\label{varest}
\Var_\infty^{z_0}(v)=\left\{\begin{array}{ll} O(|v|^{-2})\quad
&\mbox{for }\ |v|>0\\O(|v|^4\,e^{-|v|^2}) &\mbox{for }\
|v|>1\end{array}\right.\ .\end{equation}

\section{The sharp variance estimate: Proof of Theorem
\ref{sharp}}

\subsection{The codimension one case} To illustrate the basic ideas of the
argument, we begin with the proof for the case $k=1$.
 By Theorem \ref{BIPOT}, we have
\begin{equation}\label{int1s}\var\big(Z_{s^N},\phi\big)
= \int_{M}\ical^N(z)\,i\ddbar\phi(z)\;,\end{equation} where
\begin{equation}\label{int2s0}\ical^N(z)=\int_{\{z\}\times M}
Q_N(z,w)\,i\ddbar \phi(w)\;.\end{equation}

 We let
$$\Om_M=\frac 1{m!} \om^m$$ denote the volume form of $M$, and we write
\begin{equation}\label{psi}i\ddbar \phi =
\psi\,\Om_M\;,\qquad \psi \in\ccal^1_\R(M),\end{equation} so that
\begin{equation}\label{int2s}\ical^N(z)=\int_{\{z\}\times M}
Q_N(z,w)\,\psi(w)\Om_M(w)\;.\end{equation} To evaluate
$\ical^N(z_0)$ at a fixed point $z_0\in M$, we choose a normal
coordinate chart centered at $z_0$ as in \S \ref{background}, and we make
the change of variables $w=z_0+\frac v\sqrtn$. By Lemma
\ref{Qdecay} and \eqref{psi}--\eqref{int2s}, we can approximate
$\ical^N(z_0)$ by integrating \eqref{int2s} over a small ball
about $z_0$:
\begin{equation}\label{Iz}\ical^N(z_0) =   \int_{|v|\le b\sqrt{\log N}}
Q_N\left(z_0,z_0+\frac v{\sqrtn}\right)\, \psi\left(z_0+\frac
v{\sqrtn}\right)\,
\Om_M\left(z_0+\frac
v{\sqrtn}\right)+O\left(\frac 1 {N^{2m}}\right)\;,
\end{equation} where $b=\sqrt{m+2}$.

Since $\om=\frac i2 \ddbar \log a= \frac i2 \ddbar \left[|z|^2
+O(|z|^3)\right]$ in normal coordinates, we have
\begin{equation}\label{kahlerE}
\om\left(z_0+\frac v{\sqrtn}\right)= \frac
i{2}\sum\left[\de_{jk}+O\left(\frac
{|v|}{\sqrtn}\right)\right]\frac 1N dv_j\wedge d \bar v_k= \frac
i{2N}\ddbar|v|^2+O\left(\frac{|v|}{N^{3/2}}\right)\;,\end{equation}
for $|v|\le b\sqrt{\log N}$. Hence
\begin{equation}\label{OmE}\Om_M\left(z_0+\frac
v{\sqrtn}\right)=
\frac 1{m!}\left[ \frac i{2N} \ddbar |v|^2+ O\left(\frac
{|v|}{N^{3/2}}\right)\right]^m= \frac
1{N^m}\left[ \Om_E(v) +
 O\left(\sqrt\frac {\log N}{N}\right)\right]\;,\end{equation}
for $|v|\le b\sqrt{\log N}$, where $$\Om_E(v)= \frac 1{m!}\left(\frac
i2\ddbar |v|^2\right)^m= \prod_{j=1}^m\frac i2 dv_j\wedge d\bar v_j$$
denotes the Euclidean volume form. Since $\phi\in\ccal^3$ and hence
$\psi(z+\frac v{\sqrtn})=
\psi(z)+O(|v|/\sqrtn)$, we then have by Lemma \ref{FNnear} and
\eqref{Iz}--\eqref{OmE},
\begin{eqnarray}\ical^N(z_0)&=&  \frac 1 {N^m}\left[
 \int_{|v|\le b\sqrt{\log N}}  \left\{\wt G(e^{-\frac 12
|v|^2})+O(N^{-1/2+\ep})\right\}\,
\left\{\psi(z_0)+O(N^{-1/2+\ep})\right\}\right.\nonumber\\
&&\qquad \times \left\{ \Om_E(v) +
 O(N^{-1/2+\ep})\right\}\bigg]+O\left(
{N^{-m-1}}\right)\nonumber\\ &=&\frac {\psi(z_0)} {N^m}\left[
 \int_{|v|\le b\sqrt{\log N}}  \wt G(e^{-\frac 12
|v|^2})\Om_E(v) +
 O\left(N^{-1/2+\ep}\right)\right]\;.\label{U1}
\end{eqnarray}
Since $\wt G(e^{-\la})=O(e^{-2\la})$ and hence
\begin{equation}\label{U2}\int_{|v|\ge
b\sqrt{\log N}}  \wt G(e^{-\frac 12 |v|^2})\,\Om_E(v)
=O(N^{-m-1})\;,\end{equation} we can replace the integral over the
$(b\sqrt{\log N})$-ball with one over all of $\C^m$, and therefore
\begin{equation}\label{U3} \ical^N(z_0)= \frac {\psi(z_0)} {N^m}\left[
 \int_{C^m}  \wt G(e^{-\frac 12
|v|^2})\Om_E(v) +
 O(N^{-1/2+\ep})\right]\;.\end{equation}

Recalling \eqref{Gtilde}, we have
\begin{eqnarray}\label{zeta} \int_{\C^m}\wt G(e^{-\frac 12
|v|^2})\,\Om_E(v) &=& \frac 1{4\pi^2}\sum_{k=1}^\infty\int_{\C^m}
\frac{e^{-k|v|^2}}{k^2}\,\Om_E(v)\nonumber \\&=&  \frac
1{4\pi^2}\sum_{k=1}^\infty \frac {\pi^m}{k^{m+2}} \ = \ \frac
{\pi^{m-2}}4 \;\zeta(m+2) \;.\end{eqnarray}
Therefore, by
\eqref{int1s} and \eqref{U3}--\eqref{zeta},
\begin{equation}\label{U4} \var\big(Z_{s^N},\phi\big) = \frac 1{N^m}\int_M
\left[\frac{\pi^{m-2}}{4}\;\zeta(m+2) +
 O(N^{-1/2+\ep'})\right]\psi(z)^2\Om_M(z)\;.\end{equation}
The variance formula of Theorem \ref{sharp} for the case
$k=1$ follows from \eqref{psi} and \eqref{U4}.

\subsection{An explicit formula for the variance}\label{explicit}

In this section,  we give an integral formula for the variance of simultaneous
zero currents in higher codimension (Corollary \ref{BIPOTk}), which we shall
use in the next section to derive the asymptotics of Theorem \ref{sharp}.  This
integral formula is a modification of the following
formula for the variance:

\begin{theo}\label{variant}{\rm \cite[Theorem 3.13]{SZa}} Let  $1\le k\le
m$. Then for $N$ sufficiently large,
\begin{multline*} \Var\big(Z_{s_1^N,\dots,s_k^N}\big)
= \d_1\d_2\left[\sum_{j=1}^k(-1)^j{k\choose j}\dbar_1\dbar_2
Q_N\wedge\left(\d_1\dbar_1 \d_2\dbar_2 Q_{N}\right)^{j-1}\wedge \big(\E
Z_{s^N} \boxtimes \E Z_{s^N}\big)^{k-j}\right].
\end{multline*}
where  the current inside the brackets  is an $L^1$ current on $M\times M$ given by
pointwise multiplication,
$Q_N$ is given by \eqref{QN}, and  $\E Z_{s^N}$ is given by \eqref{EZ}.
Furthermore,
$\Var\big(Z_{s_1^N,\dots,s_k^N}\big)$ is an $L^1$ current on $M\times M$ if $k\le
m-1$.
\end{theo}

By an $L^1$ current on the compact manifold $M\times M$, we mean a current
whose local coefficients are $L^1$ functions.

\begin{cor}\label{variant2} Let  $1\le k\le m$.
Then for $N$ sufficiently large,
$$ \Var\big(Z_{s_1^N,\dots,s_k^N}\big)
= \d_1\dbar_1\d_2\dbar_2\left[\sum_{j=1}^k{k\choose j}
Q_N\left(-\d_1\dbar_1 \d_2\dbar_2 Q_{N}\right)^{j-1}\wedge \big(\E
Z_{s^N} \boxtimes \E Z_{s^N}\big)^{k-j}\right].$$
where  the current inside the brackets is an $L^1$ current on $M\times M$
given by pointwise multiplication.
\end{cor}

\begin{proof}     Let $$T=\sum_{j=1}^k{k\choose j}
Q_N\left(-\d_1\dbar_1 \d_2\dbar_2 Q_{N}\right)^{j-1}\wedge \big(\E
Z_{s^N} \boxtimes \E Z_{s^N}\big)^{k-j}$$ denote the expression inside the
brackets in Corollary \ref{variant2}, regarded as a $(4k-4)$-form on
$M\times M\sm \De$, where $\De=\{(z,z):z\in M\}$ denotes the
diagonal. By Lemma
\ref{DISTM2}(ii), $T=O\left(d(z,w)^{-2k+2}\right)$ and hence $T$ defines an
$L^1$ current on $M\times M$. It suffices to show that
$\dbar_1\dbar_2T$ is an $L^1$ current and thus $\dbar_1\dbar_2$ can be
moved outside the brackets in Theorem \ref{variant}.

Let $U_\ep =\{(z,w)\in M\times M: d(z,w)<\ep\}$ denote the
$\epsilon$-neighborhood of the diagonal $\De$.  For a
test form $\phi$ with $m-k+1$ $dz_j$'s, $d\bar z_j$'s, and $dw_j$'s , and
$m-k$ $d\bar w_j$'s, we  have
$$(\dbar_2T,\phi) = -\lim_{\ep\to 0} \int_{M\times M\sm U_\ep}T\wedge
\dbar_2\phi =
\lim_{\ep\to 0} \int_{M\times M\sm U_\ep}\dbar_2T\wedge \phi +
\lim_{\ep\to 0}
\int_{\d U_\ep} T\wedge \phi\;.$$  Since $T\wedge \phi =o(\ep^{-2m+1})$
on $\d U_\ep$, the boundary integral goes to 0.   By Lemma \ref{DISTM2}(ii)
and the fact that $Q_N\in\ccal^1(M\times M)$, the pointwise-defined form
$\dbar_2T$ is also $O(d(z,w)^{-2k+2})$ and thus $\dbar_2T$ is
an $L^1$ current on $M\times M$ given by pointwise differentiation.

Repeating the same argument with $T$ replaced by $\dbar_2T$ and using part (i) of
Lemma \ref{DISTM2} as well as part (ii), we then conclude that $\dbar_1\dbar_2T$ is an
$L^1$ current.\end{proof}

Corollary \ref{variant2} can also be shown directly,
using the argument in the proof of Theorem \ref{variant}
in \cite{SZa}.

\begin{cor}
\label{BIPOTk} The variance in Theorem \ref{sharp} is given by:
\begin{multline*} \var\big([Z_{s_1^N,\dots,s_k^N}],\phi\big)\\=
 \sum_{j=1}^k{k\choose j}\int_{M\times M}Q_N\left(-\d_1\dbar_1
\d_2\dbar_2
 Q_{N}\right)^{j-1}\wedge
\big(\E Z_{s^N} \boxtimes \E Z_{s^N}\big)^{k-j}\wedge (
\ddbar\phi\boxtimes \ddbar\phi),
\end{multline*}
where  the integrands  are in $L^1(M\times M)$.
\end{cor}

\subsection{Higher codimensions}\label{higher} Recalling
\eqref{indepE}, we  write the  formula of Corollary
\ref{BIPOTk} as follows:
\begin{eqnarray}\label{Vj0} \var\big([Z_{s_1^N,\dots,s_k^N}],\phi\big) &=&
\sum_{j=1}^k {k\choose j}\, V_j^N(\phi)\;,\\V_j^N(\phi)&=&
\left(\frac N\pi \right)^{2k-2j} \int_{M\times M} Q_N(z,w)\,\big[-
\d_1\dbar_1 \d_2\dbar_2 Q_{N}(z,w)\big]^{j-1}\nonumber\\&&
\quad\wedge\left[\om(z)^{k-j}+O\left(\frac
1N\right)\right]\wedge\left[\om(w)^{k-j}+O\left(\frac
1N\right)\right]\wedge i\ddbar\phi(z)\wedge i\ddbar\phi(w)
\nonumber\\&=&\left(\frac N\pi \right)^{2k-2j} \int_M
\ical^N_j\wedge\left[\om^{k-j} \wedge  i\ddbar\phi+O\left(\frac
1N\right)\right]\,,\label{Vj}\end{eqnarray} where
\begin{eqnarray}\label{Ij}\ical^N_j(z) &=& \int _{\{z\}\times
M}Q_N(z,w)\,\big[-
\d_1\dbar_1 \d_2\dbar_2
 Q_{N}(z,w)\big]^{j-1}\wedge\left[
\om(w)^{k-j} \wedge i\ddbar\phi(w)+O\left(\frac
1N\right)\right]\nonumber\\ &&\qquad\in
T^{*j-1,j-1}_z(M)\;.\end{eqnarray} The integrand in \eqref{Ij} is
regarded as an $(m,m)$-form (in the $w$ variable) with values in
$T^{*j-1,j-1}_z(M)$.

Fix a point $z_0\in M$, and let $1\le j\le k\le m$.  To evaluate
$\ical^N_j(z_0)$, we write
\begin{eqnarray*}\om^{k-j}\wedge i\ddbar\phi &=& \sum \psi_{JK}(w)
dw^J\wedge d\bar w^K\\
&=&\frac{1}{N^{m-j+1}}\sum \psi_{JK}\left(z_0+\frac v\sqrtn\right)
dv^J\wedge d\bar v^K\;,\qquad |J|=|K|=m-j+1\;.\end{eqnarray*}

By Lemma \ref{Qdecay}, we can replace integration over $M$ in
\eqref{Ij} with integration over the small ball of radius
$b\sqrt{\log N/N}$, with $b=\sqrt{m+4}$, to obtain:
\begin{eqnarray*} \ical^N_j(z_0)
&=&N^{-m+j-1}\int_{|v|\le b\sqrt{\log N}}  Q_N\left(z_0,z_0+\frac
v \sqrtn\right)\,\left[- \d_1\dbar_1 \d_2\dbar_2
Q_{N}\left(z_0,z_0+\frac v \sqrtn\right)\right]^{j-1}\\&&
\qquad\wedge \sum\left[ \psi_{JK}\left(z_0+\frac v\sqrtn\right)
+O\left(\frac 1N\right)\right]dv^J\wedge d\bar v^K
+O\left(N^{-m-1}\right)\;.\end{eqnarray*} By Lemma \ref{DISTM2},
the above integrand is $L^1$, and hence by Lemma \ref{d4Qas},
\begin{eqnarray} \label{Ijapprox}\ical^N_j(z_0) &=&
N^{2j-2-m}\left[ \int_{|v|\le b\sqrt{\log N}} F(\half
|v|^2)\Big\{\Var^{z_0}_\infty(v)\Big\}^{j-1}\sum
\psi_{JK}(z_0) dv^J\wedge d\bar v^K +O(N^{-\frac 12+\ep})\right]\nonumber\\
&= & N^{2j-2-m} \left[\sum \psi_{JK}\left(z_0\right)
\int_{v\in\C^m}\!\! F(\half |v|^2) \Big\{
\Var^{z_0}_\infty(v)\Big\}^{j-1} dv^J\wedge d\bar v^K+O(N^{-\frac
12+\ep})\right].\end{eqnarray} Here, we replaced the integral over
the $(b\sqrt{\log N})$-ball with one over all of $\C^m$, since  by
\eqref{varest} we have $F(\half |v|^2) \{
\Var^{z_0}_\infty(v)\}^{j-1} =O(e^{-|v|^2})$ for $|v|>1$, and
hence
$$\int_{|v|>b\sqrt{\log N}}F(\half |v|^2) \Big\{
\Var^{z_0}_\infty(v)\Big\}^{j-1} dv^J\wedge d\bar
v^K=O(N^{-4m})\;.$$ It follows from \eqref{Vj} and
\eqref{Ijapprox} that
\begin{equation}\label{Bj} V_j^N(\phi)=
N^{2k-m-2}\left[\int_M\sum\ B^j_{JKAB} \psi_{JK}\bar \psi_{AB}
\,\Om_M+O(N^{-\frac 12+\ep})\right]\;,\end{equation} where
$B^j=\{B^j_{JKAB}\}$ is a universal Hermitian form on
$T^{*m-j+1,m-j+1}(M)$. Theorem \ref{sharp} then follows from
\eqref{Vj0} and \eqref{Bj} with $$B_{mk}(\al,\al)=\sum_{j=1}^k
B^j(\om^{k-j}\wedge \al,\, \om^{k-j}\wedge \al)\;.$$\qed

\section{Asymptotic normality: Proof of Theorem
\ref{AN}}\label{s-normality}

The proof of Theorem
\ref{AN} is an application of  Propositions \ref{DPdecay}--\ref{better}
to a general result of Sodin-Tsirelson \cite{ST} on asymptotic normality
of nonlinear functionals of Gaussian processes. Following \cite{ST}, we
define a {\it normalized complex Gaussian process\/} to be a
complex-valued random function $w(t)$ on a measure space $(T,
\mu)$ of the form
$$w(t)=\sum c_j  g_j(t)\;,$$ where the $c_j$ are i.i.d.\ complex
Gaussian random variables (of mean 0, variance 1), and the $g_j$ are
(fixed) complex-valued measurable functions such that
$$\sum|g_j(t)|^2=1\quad
\mbox{for all }\ t\in T.$$
We let $w^1,w^2,w^3,\dots$ be a sequence of normalized complex
Gaussian processes on a finite measure space $(T, \mu)$. Let $f(r)
\in L^2(\R^+, e^{-r^2/2} rdr)$ and let
$\psi:  T
\to
\R$ be bounded measurable. We write
$$Z_N^{\psi}(w^N) = \int_T f(|w^N(t)|) \psi(t) d\mu(t).$$

\begin{theo}\label{ST} \cite[Theorem 2.2]{ST} Let $\rho_N(s, t)$ be the
covariance functions for the Gaussian processes $w^N(t)$. Suppose that
\begin{enumerate}

\item[i)]  \quad $\displaystyle\liminf_{N \to \infty} \frac{\int_T
\int_{T} |\rho_N(s, t)|^{2 \alpha} \psi(s) \psi(t) d\mu(s)
d\mu(t)}{\sup_{ s \in T} \int_T |\rho_N(s, t)| d\mu(t)} > 0\;,$\\ for
$\al=1$ if  $f$ is monotonically increasing, or for all $\al\in\Z^+$
otherwise;

\item[ii)]  \quad $\displaystyle\lim_{N \to \infty}\; \sup_{s \in T}
\int_T |\rho_N(s, t)| d\mu(t) = 0.$

\end{enumerate}
Then the distributions of the random variables
$$ \frac{Z_N^{\psi} - \E Z_N^{\psi}}{\sqrt{\var(Z_N^{\psi})}}$$
converge weakly to
$\ncal(0, 1)$ as $ N \to \infty$.
\end{theo}

We apply this result with $f(r)=\log r$ and $(T, \mu) = (M,
\Om_M)$. To define our normalized Gaussian processes $w^N$ on $M$, we
choose a measurable section $\sigma_L:M\to L$ of $L$ with
$\|\sigma_L(z)\|_h=1$ for all $z\in M$, and we let
$$S_j^N=F^N_j\sigma_L^{\otimes N}\,,\quad j=1,\dots,d_N,$$ be an
orthonormal basis for $H^0(M, L^N)$ with respect to its Hermitian
Gaussian measure, for each $N\in\Z^+$.  We then  let
$$g^N_j(z) : = \frac{F^N_j(z)}{\sqrt{\Pi_N(z,z)}}\,,\quad
j=1,\dots,d_N\,.$$
Since $|F^N_j|=\|S^N_j\|_{h^N}$, it follows that $w^N= \sum c_j  g^N_j$
defines a normalized complex
Gaussian process, for each  $N\in\Z^+$ (where the $c_j$ are i.i.d.\
standard complex Gaussian random variables). In fact,
$$|w^N(z)| = \frac{\|s^N(z)\|_{h^N}}{\sqrt{\Pi_N(z,z)}}, $$
where $$s^N= \sqrt{\Pi_N(z,z)}\,w^N\sigma_L^{\otimes N}= \sum c_j S^N_j$$
is a random holomorphic section in $H^0(M, L^N)$. The covariance functions
$\rho_N(z,w)$ for these Gaussian processes satisfy
$$|\rho_N(z,w)|=P_N(z,w)\;.$$

We now let $\phi$ be a fixed  $\ccal^3$ real
$(m-1,m-1)$-form on $M$ and we write
$$\frac i\pi\ddbar\phi =\psi\,
\Omega_M\;.$$ Then $\psi \in \ccal^1$, and
$$Z_N^{\psi}(w^N) = \int_M\left( \log \|s^N(z)\|_{h^N} - \log
\sqrt{\Pi_N(z,z)} \right)\frac i\pi \ddbar \phi(z) = \big(
Z_{s^N},\phi\big) + k_N\,,$$
where the $k_N$ are constants (depending on $L\to M$ and $\phi$, but
independent of the random sections $s^N$). Hence $Z_N^{\psi}(w^N)$ has the
same variance as the smooth linear statistic $\big(Z_{s^N},\phi\big)$, and
it suffices by Theorem \ref{ST} to check that the covariance function
satisfies conditions (i)--(ii) of the theorem. We start with (ii): by
Proposition
\ref{DPdecay},
$$\lim_{N \to \infty}\; \sup_{z \in M} \int_{d(z,w) > b\sqrt{\frac
{\log N}{N}}} P_N(z,w) \,\Omega_M(w)= 0. $$ On the other hand,
since $P_N(z,w) \leq 1$, it is obvious that the same limit holds
for $d(z,w) \leq b\sqrt{\frac{ \log N}{N}}$,  verifying (ii).

To check (i), we again break up the integral into the near diagonal
$d(z,w) \leq b\sqrt{\frac{ \log N}{N}}$ and the off-diagonal $d(z,w) >
b\sqrt{\frac{ \log N}{N}}$. As before, the integrals over the off-diagonal
set tend to zero rapidly and can be ignored in both the numerator and
denominator.

On the near diagonal, we  replace $P_N$ by its asymptotics in
Proposition \ref{better}. The asymptotic formula for $P_N$ has a
universal leading term independent of $z$ and has uniform remainder, so
condition (i) (with
$\al=1$) becomes
$$\liminf_{N \to \infty} \frac{\int_M \Omega_M(z)\int_{|u| < b\sqrt{\log N}
} e^{- |u|^2}[1 + R_N(u)]^{2 } \psi(z +
\frac{u}{\sqrt{N}}) \psi(z)\, du}{ \int_{|u| <b
\sqrt{\log N} } e^{-\half |u|^2}[1 + R_N(u)]\, du} >
0. $$ Since $\psi \in \ccal^1$, the ratio clearly tends to $2^{-m}\int_M
\psi(z)^2 \,\Omega_M > 0$, which verifies (i) and completes the
proof of Theorem \ref{AN}.\qed

\section{Open problems on smooth and counting
statistics}\label{OPEN}

In this section, we present a number of open problems on
smooth and discontinuous linear statistics of
zeros.

\begin{enumerate}

\item Asymptotic normality of the smooth linear statistics
$(Z_{s_1^N,\dots,s_k^N},\phi)$ has only been proved in codimension
one, i.e.\ when $k = 1$. But these random variables are likely
to be asymptotically normal for all dimensions $m$ and
codimensions $k$. It would be interesting to prove (or disprove)
this statement, in particular for $k = m$.

\item To our knowledge, no results are known  to date
regarding  the asymptotic normality of the  counting statistics
$\ncal_N^U(s_1^N,\dots,s_m^N)$, even when $m = 1$.  This is
analogous to, but presumably harder than, the smooth linear
statistic when $k = m$.

\item In \cite{Zh}, Qi Zhong obtained surprising results on the
expected value of the `energy' random variable
\begin{equation} \label{ENERGY} \ecal_G (s^N) = \sum_{i \not = j}
G(a_i, a_j), \quad  Z_{s^N}= \{a_1,\dots,a_N\}\,, \end{equation}
summing the values of the Green's function $G$ over pairs of
distinct zeros of a random polynomial or section  $s^N \in H^0(M, L^N)$
of the $N$-th power of a positive line bundle $(L, h)$
over a compact Riemann surface $M$. (The Green's function is normalized to
equal $+ \infty$ on the diagonal.)
Zhong proved  that when $G$ is the Green's function for the
Riemannian metric induced by the curvature of $h$, the expected
energy has the asymptotics $\E \ecal_G (p_N) \sim -\frac{1}{4 \pi}
N \log N$. It  is known (N. Elkies) that $-\frac{1}{4 \pi} N \log
N$ is also the asymptotic minimum for the energy sum (\ref{ENERGY}).  The
energy is a partially smooth linear statistic \begin{equation} \label{G} (G,
Z_{s^N} \otimes Z_{s^N} - \Delta_{Z_{s^N}} )\end{equation}
on $M \times M$, where $\Delta_{Z_{s^N}}$ represents the diagonal
terms of
 $Z_{s^N}\otimes Z_{s^N}$. The statistic (\ref{G})  is not  smooth since $G$ has a
logarithmic singularity and since  we subtracted the diagonal
current. It is a random variable of one section $s_N$ in
dimension one, but $G$ is a function on $M \times M$, so the variance
of $\ecal_G$ involves the rather complicated  four-point correlation function
of $Z_s$ rather than the pair correlation, for which Theorem \ref{BIPOT}
gives a useful formula. It would be interesting to investigate the variance
of the energy $\ecal_G$. It seems that it should tend to zero with $N$  since
the Gaussian measure is  concentrated on `polynomials' whose zeros are
asymptotic minimizers.

\item The expected distribution of zeros can have quite disparate
asymptotics when the ensembles are given Gaussian measures induced from
 inner products on the space of polynomials (or sections) which use non-smooth volume
 forms or non-positively curved line bundle metrics. For instance, in the case
 where the measure is supported on an analytic plane domain $\Omega \subset \C$
 or on its boundary, it was shown in \cite{SZeq} that the  the expected distributions of random zeros
 of random polynomials of degree $N$ tend to the equilibrium measure of $\Omega$. 
 This result was generalized to higher dimensions and more general  metrics and measures in a sequence
 of papers \cite{Bl,BS,Be} in which it is shown that  the expected distribution  of zeros tends to an equilibrium measure adapted
 to the measure and metric. In \cite{Sh},
an upper bound was given for the variances of the smooth linear statistics
when the inner products are defined by arbitrary measures (and also for more
general sequences of ensembles of increasing degrees). This bound is
sufficient to prove that sequences of random zeros in these ensembles almost
surely converge to their equilibrium measure, although the bound is not
always sharp.  The \szego kernels for the inner products in 
\cite{SZeq,Bl,BS,Be} are quite  different from \szego kernels
for  positive line bundles in this article, 
 and so the asymptotics of the variances might be quite different. It would be interesting
to determine them.

\item   Results on expected values for the analogous SO$(2)$ and SO$(m+1)$
ensembles of random real polynomials in one or several variables were
given by \cite{BD, EK,  Ro, SS}, and a (global) variance result for real
zeros was given in \cite{Ws}.  Maslova \cite{Ma} proved the asymptotic
normality of  the number of real zeros for the Kac ensemble \cite{Kac}
of random real polynomials on
$\R^1$  (as well as for some non-Gaussian ensembles). But as far as
we are aware, asymptotic normality for numbers of real zeros in the
SO$(m+1)$ ensemble has not been investigated.

\end{enumerate}

\end{document}